\documentclass{amsart} 

\usepackage{amsmath, amsthm}
\usepackage{amssymb, latexsym}
\usepackage{amsfonts} 
\usepackage[dvipsnames]{xcolor} 
 
\usepackage{threeparttable} 

    \usepackage{lineno} 
\newtheorem{Theorem}[equation]{Theorem}

\newtheorem{Proposition}[equation]{Proposition}
\newtheorem{proposition}[equation]{Proposition} 
 
\newtheorem{lemma}[equation]{Lemma} 

\newtheorem{Corollary}[equation]{Corollary}
\newtheorem{corollary}[equation]{Corollary}

\newtheorem{Examples}[equation]{Examples}

\newtheorem{Remark}[equation]{Remark}

\theoremstyle{definition}

\newcommand{\exsref}[1]{\rm Ex\-am\-ples \ref{#1}}
\newcommand{\thmref}[1]{\rm Theorem~\ref{#1}}

\newcommand{\secref}[1]{\rm Sec\-tion \ref{#1}}
\newcommand{\lemref}[1]{\rm Lem\-ma \ref{#1}}
\newcommand{\propref}[1]{\rm Proposition~\ref{#1}}
\newcommand{\corref}[1]{\rm Corollary~\ref{#1}}

\newcommand{\remref}[1]{Re\-mark \ref{#1}}

\newcommand{\beql}[1]{\begin{equation}\label{#1}}
\newcommand{\eeq} {\end{equation}}

    \font\Aaa=msam10
\newcommand{\nor}{\mbox{ $\vartriangleleft$ }}
\newcommand{\ron}{\mbox{ $\vartriangleright$ }}

\font\Bbb=msbm10

\newcommand\Z{\hbox{\Bbb Z}}

\newcommand\F{\hbox{\Bbb F}}

\numberwithin{equation}{section}

\DeclareMathOperator{\Aut}{Aut}

\newcommand\G{\Gamma }
\newcommand\D{{ \Delta }}

\DeclareMathOperator{\AGL}{{\rm A\Gamma L}}

          \def\AgL{\mbox{\rm A$\Gamma$L}}
        \def\qed{\hbox{~~\Aaa\char'003}}
        \def\AGL{{\rm AGL}}

        \def\<{{\langle}}
        \def\>{{\rangle}}

 \def\a{\alpha}

\def\div{ \kern-.5pt\hbox{\big |} }
\def\ndiv{ {\not\kern-.5pt\hbox{\big |}\,} }
\def\ndivv{ {\not\kern+1.5pt\hbox{$\mid$}\,} }

\def\Aut{{\rm Aut}}

\def\col{\colon\!}

\def\B{^2\kern-.8pt B}
\def\G{^2\kern-.8pt G}
\def\EH{^2\kern-.8pt\hat  E}
\def\E{^2\kern-.8pt E}
\def\D{^3\kern-1pt D}
\def\FF{^2\kern-.8pt F}

\newdimen\refcodesize
\newbox\seriesbox
\def\proof{\smallskip\noindent {\bf Proof.~}}

 \font\Bbb=msbm10

 \font\Bbb=msbm10

\def\Z{\hbox{\Bbb Z}}

\def\F{\hbox{\Bbb F}}

\def\phi{\varphi}

\def\col{\colon\!}

\DeclareRobustCommand{\SkipTocEntry}[4]{}
 
\usepackage{lineno}
 \numberwithin{equation}{section}

\begin{document}

\title[Soft planes and groups]
{Soft planes and groups
 }

       \author{William  M. Kantor
      }
      \address{Brook House, Brookline, MA 02445}
  \email{kantor@uoregon.edu}

\begin{abstract}

Finite projective planes are constructed using groups that
 satisfy simple-looking conditions.   
 The resulting projective planes 
include many  known planes and possibly new   ones, and
 are precisely those having 
 a collineation group   fixing a flag $(\infty,L_ \infty )$  and transitive on the  flags 
 $ (w,W )$ with $w\notin L_ \infty $ and $\infty\notin W$.
   
 \vspace{-18pt}  
 \end{abstract}
\maketitle
\section{Introduction
   }
 \label{Introduction} 
  
Finite projective planes are constructed using groups.
This does not even require  knowing what a finite projective plane is:  
what is needed is the construction of groups that
 satisfy simple-looking conditions
(Theorem \ref{main}).

The  planes constructed are {\em soft planes}:  ``sort-of-flag-transitive planes'',
meaning that there is   a collineation group  ({\em a soft group})
 fixing a flag $(\infty,L_ \infty )$  and transitive on the  flags 
 $( w,W )$ with $w\notin L_ \infty $ and $\infty\notin W$.  
There are   known families of examples of such planes.
The initial goal was to obtain examples not of prime power order
(contradicting a standard conjecture \cite[p.~294]{Pi}, \cite[p.~25]{Ry}), 
but~none has yet been found.  A secondary goal  was to construct examples of prime power order that are neither translation planes nor their duals, but that was also not successful.%

\secref{Construction of planes} proves the relationship between soft planes and groups
using undergraduate algebra, while 
\secref{properties} contains elementary observations concerning the planes and the groups.
\secref{summary} summarizes the purely group-theoretic restrictions presently known,
such as those involving  normal structure.
\secref{Additional groups} contains additional   soft groups for some of  the
known soft planes.
This paper has an old-fashioned  point of view:    its 
methods are in the 57-year-old book \cite{De}.  
 
\smallskip  \smallskip
{\noindent\bf Background.}  
Higman and McLaughlin \cite{HM} constructed a flag-transi\-tive geometry using the cosets of two subgroups $A$ and $B$ of a group $G$, with a ``point'' $Ax$ ``on'' a ``line'' $By$ iff  $Ax\cap By\ne \emptyset$.\footnote{In a much more general context the geometric  use of cosets and their intersections    seems to have originated in \cite[p.~272]{Ti}.}
  One of their goals was 
to study finite flag-transitive  projective planes.
The present paper arose from the realization that some ideas in \cite{HM} could 
  be used while assuming less transitivity.
The ``new'' idea is in Theorems~\ref{main}  and \ref{converse},  based on  \cite[Lemma~4]{HM}.
In the 1960s this idea might have seemed novel, but now it appears straightforward.
It seems surprising that   \thmref{main}  has not previously been observed.  

I am indebted to Theobald Grundh\"ofer for  helpful information about the geometry of points and lines of the underlying soft plane that are on neither $\infty$ nor $L_\infty$, and
on which a soft group acts flag-transitively.   Historical comments   are in \cite[pp.~228-229]{BHKT}; the geometry was placed in a larger context in 
\cite[pp.~315-316]{De}.  The only reference  mentioning the collineation group  of such a geometry  appears to be \cite{Cr} for  \exsref{Heisenberg}.
 
 \section{Construction of planes}

\label{Construction of planes}

This paper   concerns the  following result (and its converse, \thmref{converse}).~The 
proof  involves little more than the definition of a projective plane and 
elementary properties of  cosets.
\Theorem
\label{main}
There is a projective plane $\pi$ of order $n>1$   if there is a group $G$ 
having subgroups $A ,$  $B$  and $M$ such that 
 \vspace{-2pt}
 \begin{enumerate} 
 \item[\rm(1)]     $|A|=|B|=|M|=  nk  $   and  $|G|=n^3k$ for $k=|A\cap B| ,$    
  \item[\rm(2)]  $AM$  and $BM$ are subgroups of  order  $n^2k,$  
  \item[\rm(3)]    $G= AMB ,$ 
   and 
\item[\rm(4)]$AB\cap BA=A\cup B$.
 \vspace{-2pt}
        \end{enumerate} 
Moreover$,$ $\pi$   is a  soft plane  with  $G$ inducing a soft group.
  \smallskip\smallskip
       \rm
       
       Here  $AM:=\{a m\mid a\in A, m\in M\}$.     
      
      \Remark\rm
The intriguing aspect of this result is that there is no 
obvious arithmetic reason to expect $n$ to be a prime power.

Even though a question about  some projective planes is equivalent to one about some groups, the latter setting ought to provide many more tools, involving (normal)
subgroups, quotient groups, group extensions, etc.  However, as indicated earlier, these additional facilities have  not yet  led to any new planes.

The principal obstacle is the  condition (4).

      \Remark\rm\label{faithful}
In the theorem $G$ might not act faithfully on $\pi$.  
{\em After the proof we will always assume that $G$ 
acts faithfully as
a group of collineations of $\pi$.}

It seems likely that $G$ is then solvable (cf. \propref{M elations}(iii)).
 
       \Examples[\rm Heisenberg groups]\rm
       \label{Heisenberg}  \hspace{-4pt}
  \hspace{-2pt} Using entries from $\F_n$ when $n$ is a prime~power$,$%
\smallskip\smallskip  

$G:=   \begin{pmatrix} 
      1 & *& * \\
      0 & 1& * \\
      0 & 0& 1 \\
   \end{pmatrix}\! ,  \
   A:=
   \begin{pmatrix} 
      1 & *& 0 \\
      0 & 1& 0 \\
      0 & 0& 1 \\
   \end{pmatrix}\! ,  \
  M:=
   \begin{pmatrix} 
      1 & 0& * \\
      0 & 1& 0 \\
      0 & 0& 1 \\
   \end{pmatrix}\!   $ \   and
 $ B:=
   \begin{pmatrix} 
      1 & 0& 0 \\
      0 & 1& * \\
      0 & 0& 1 \\
   \end{pmatrix}
   $
   
   \smallskip\smallskip
   \noindent behave as in the theorem; they produce the desarguesian plane of order $n$.
     More  general  
 Heisenberg groups  arise   by replacing $\F_n$ by any finite 
 {\em semifield} (nonassociative division algebra) \cite{Cr,Hi};
 these groups have class 2.
     
 \Examples\rm
 \label{list}
  Up to duality (i.\,e., interchanging $A$ and $B$) the  
  {\em known} planes    arising via (1,2,3,4) are the following translation planes:%
\begin{enumerate}
 \vspace{-2pt}
\item[\rm(i)] semifield planes \cite{Cr,Hi} (\exsref{Heisenberg}),  
  
  \item[\rm(ii)]    likeable planes  with $n=q^2$  
  when  $q>2$  is a  prime power with ${q\equiv 2}$~(mod~3) or 
  when $q>5$ is a power of 5  \cite{Ka}  (cf.~\cite{Be,Wa}),  

\item[\rm(iii)] a plane due to Sherk     \cite{Sh,NRS}   with  $n=3^3$,
\item[\rm(iv)]  L\"uneburg-Tits planes  \cite{Lu}  with $n=q^2$ for $q=2^{2e+1}>4$, 

\item[\rm(v)] a plane  due to Jha-Johnson  \cite{JJ}    with  $n=8^2$
(cf.  \cite[p.~13]{BJJM}),
 and
\item[\rm(vi)] a plane due to Biliotti-Menichetti      \cite{BM}  
 with  $n=8^2$.
  \vspace{-2pt}
\end{enumerate}
Each plane    (ii-vi)  was  found at least 40 years~ago without the use of a computer.
While all of these examples were found using soft groups with $k=1$ (cf.    \remref{further}),
many of these planes have soft groups with $k>1$  (cf. \remref{solvable}).
   
  \Remark\rm 
  \label{further}
We briefly describe further properties   of  some of the   examples in  (i-vi) 
{\em  for which $M\nor G,$  $k=1$ and $G$ is a soft $p$-group of order $n^3$.}

 In (i)  $M=Z(G)$  and  there are many possible nonisomorphic groups $G$
 \cite{Hi},  \cite[Sec.~5.3]{De}. 
 For each of them  $G/ M$ is elementary abelian.
  
  In (ii) $A$ is abelian: elementary abelian when $n$ is odd but   homocyclic of exponent 4 when $n$ is even.

 In (iii)   $A \cong \Z_3\times\Z_9$.  
 
    In (iv)  $A$  is a Sylow 2-subgroup of a Suzuki  group  ${\rm Sz}(q)$
   and $BM$ is the natural 4-dimensional module for ${\rm Sz}(q)$.
 
     In (v,vi)   $A$ is nonabelian  \cite[pp.~12-13]{BJJM}.  
  
  The planes  (ii-vi) were   obtained by requiring   $BM$ 
  to be an elementary abelian normal subgroup of $  G$ (but 
  $AM \ntrianglelefteq G$).  Moreover,
  $G \lesssim \AgL(4 ,q)$   in (ii,iv,v,vi) 
     or $\AGL(6,3)$~in~(iii),    
$G$ is a $p$-group of   nilpotence
       class   4 or 5,  and $Z(G) = C_M(A)$ has order
  $\sqrt n$  in (ii,iv,v,vi) or      3  in (iii).
  
While $G$ uniquely determines $\pi$, \secref{Additional groups}  
contains additional   examples of groups $G$ that
determine some of the above planes but have different normal structures:
  examples having 
   $M\ntrianglelefteq G$,
$AM\ntrianglelefteq  G$ and $BM\ntrianglelefteq  G,$ and, what seems more interesting, examples of order $n^3$  having  $M\nor G$, 
 $AM\ntrianglelefteq  G$ and $BM\ntrianglelefteq  G$.
          \vspace{3pt}
          
{\noindent \bf  Proof of  Theorem \ref{main}.}
 The plane $\pi$ can be described as follows:
  
 \hspace{24pt}  distinguished line $L_\infty$ and distinguished point $\infty $ on it
  
  \hspace{24pt}   other points:  cosets $Ax$ and cosets $BMx$ \ \hspace{-0.9pt}for  $x \in G$

 \hspace{24pt}  other lines:  \ \hspace{-.5pt} cosets $By$ and cosets $AMy$ \ for $y\in G$
  
 \vspace{-2pt} 
 \noindent
 \begin{minipage}{.57\linewidth}
 \begin{equation}
 \label{table}
\begin{array}{lllllll}
 Ax  &\hspace {-6pt}is ~on &\hspace {-8pt} By  & \hspace {-8pt} \iff Ax\cap B y \ne \emptyset 
 \\ 
 Ax  &\hspace {-6pt} is ~on &\hspace {-8pt}AMy &\hspace {-8pt}\iff Ax\subset AMy
 \\
 By & \hspace {-6pt}is ~on  &\hspace {-8pt}BMx     &\hspace {-8pt}\iff By\subset BMx 
\end{array} 
 \end{equation}
  
  \vspace {-3pt} 
\em
 \hspace {34pt} $BM x$   is on   $L_\infty$
\ and \
$AM y$    is on   $\infty$.
  \rm
  \end{minipage}

  \rm
    \vspace {3pt}  

Thus, 
{\em a point-coset is on a line-coset if and only if the cosets meet.}

That  
 $Ax $ is  on  $ By    \iff Ax\cap B y\ne \emptyset $
 is taken from \cite[Lemma~4]{HM}.

Here and later    $A,B,AM$ and $BM$ have different lives as   
subgroups of $G$  and  as points or lines of  $\pi$, which is occasionally awkward.
 
 {\em
Right multiplications by elements of $G$ induce automorphisms of $\pi$.}  

The structure of $G$ is not involved in the proof, which only uses 
 the most basic properties of cosets: {\em the proof is elementary.}

It is easy to use (1) and (2) in the theorem to check that there are $n^2+n+1$ points and lines, 
with $n+1$ points per line and $n+1$ lines per point.   
We   check that every two 
points are on one and only one line  in a few simple steps involving \eqref{table}.%
 \vspace{1pt} 

(i)   
First note that    $AM\cap B= A\cap B$  (since 
  $AM\cap B\ge  A\cap B$ and, by   (1,2,3),   
 $|B|/|A\cap B|=n=|G|/|AM|=|B|/|AM\cap B|$).   Similarly,  $BM\cap A= A\cap B$.
 \vspace{1pt}

(ii) {\em Distinct points $A$ and $Am,$ $m\in M,$ are on $AM$ and no other line.}
(If $A$ and $Am$ are on $Ba$ with  $a\in  A$  then  $a'm=b'a$ with $a'\in A, b'\in B$.  Then 
$b'\in  AM\cap B=A\cap B$ by (i), so $a'm\in A$ and $Am=A$.)
 \vspace{1pt}

 (iii)  By (4),    {\em two points not on $L_\infty$ are on at most one line}  
\cite[Lemma~4]{HM}.
  (For, we may assume that we have two points $A $ and $Ab $  of  $B$
 with $ b \in B\backslash  A $.  They are not on $AM$ since 
  $b\in AM\cap B\implies b\in A\cap B$ by (i). 
 If they are on a line $Bz $ then 
 $z=b_1a_1 =b_2a_2b $ for some $b_1 ,b_2 \in B$ and  $ a_1,a_2\in A$.  
 Now $(b_2^{-1}b_1)a_1=a_2b\in $ $BA\cap AB=A\cup B$ by (4), so
  $  b \notin A \implies (b_2^{-1}b_1)a_1\in B$ and   $Bz=   B a_1  =B.$)%
  \vspace{1pt}
 
(iv) {\em The points  $A$  and $\infty$ are on a unique line} ($A\subset AMy \iff AMy= AM$).
 \vspace{1pt}

 {\em The points $A$ and $BMx$ are on a unique line.}   (By (3),  $BMx=BMa$ with~$a\in A$,
so  $A$ and  $BMx$  are on~$Ba$. 
By \eqref{table}, a line on $A$ and $BM x$ must look like $Ba' \subset  BM x=BMa$  with  $ a'\in A$, 
 so $a'a^{-1}\in BM\cap A=A\cap B$ by (i), and $Ba'=Ba$.)%
  \vspace{1pt}
 
 (v) We know that $A$ is on a unique line with each of the $(n-1)+1 $ points  $\ne A$ on $AM$ and with each of  the
 $n((n-1)+1 )$ points $\ne A$  on the  $n$  lines 
  $Ba,$ $ a\in A$,
 and hence with each of   $n^2+n$   points  $\ne A$.~Thus,  {\em$\pi$ is a projective plane of order~$n$.}%
  \vspace{1pt}
 
  (vi) 
 Finally, {\em$G$ acts as a soft group for $\pi$}:  it is transitive on the flags 
 $(Ax,By)$ with $x,y\in G$.
 (For $(Ax,By)$  to be a flag we must have $ax=by$ for some $a\in A,b\in B$,
 so $(A,B)ax=(Ax,By)$.)   \qed
   
\Remark \label{ideals}\rm 
 In this paper $\mathfrak  I$ denotes the  set of points on the line of $\pi$ on  both  $ A$ 
 and $ \infty$  $($an {\em ideal line}$) ,$  and  
   $\mathfrak  i$  denotes the set of lines   on  the point  on both  $B$ and $L_{ \infty}$
 $($an {\em ideal point}$) $.   
 $($We will also view $ \mathfrak  I $ as a line of $\pi$  and $\mathfrak  i$  as a point.$)$
 If $G_{\mathfrak  I} $ and $G_{\mathfrak  i} $  denote  the  stabilizers
  of ${\mathfrak  I} $    and ${\mathfrak  i}, $ respectively$,$  then
  \vspace{-2pt}
 \begin{itemize}
  \item [(i)]   ${\mathfrak  I}  = \{Am\mid m\in M\}  \dot\cup \{ \infty \} $  and
  ${\mathfrak  i}  = \{Bm\mid m\in M\}  \dot\cup \{L_ \infty \},$ 
   \item [(ii)]   $AM=G_{\mathfrak  I} $ and  $BM=G_{\mathfrak  i},$  
    \item [(iii)] $M=G_{\mathfrak  I\mathfrak  i} ,$ and
       \item [(iv)] 
   $BM$ is transitive on the points not on $L_ \infty$  and   $AM$ is transitive on the 
   lines not   on  $\infty.$
  \end{itemize} \rm
 
 \proof   We will use \eqref{table}.
 
 (i) ${\mathfrak  I}$ is the set of points on the line $AM$.
  \vspace{1pt}
  
 (ii) $G_{\mathfrak  I} =G_{AM} =AM$.
  \vspace{1pt}
 
 (iii)  By \eqref{table}, \thmref{main}(3) and (ii), $G_{\mathfrak  I}$ is 
  transitive on $L_ \infty \backslash \{\infty\}.$
  By (ii) and \thmref{main}(1,2),
 $|G_{\mathfrak  I \mathfrak i} |= |AM|/n=n^2k/n=|M|$  
 and $G_{\mathfrak  I \mathfrak i}\ge M$.

 (iv)  By \thmref{main}(2,3), $G=A(BM) =B(AM) $.
  \qed
 
\Corollary
\label{more conditions} 
The subgroups 
$A,$ $B$ and $M$   in {\thmref{main}} also satisfy 
\vspace{-2pt}
\begin{itemize}
 \item [(i)] $A\cap M=A\cap B=B\cap M,$  
  \item [(ii)]  $ AM\cap B=A\cap B =B M\cap A,$   
 \item [(iii)]   $AM\cap BM= M,$    
 \item [(iv)]  $AM  =G\backslash \big(A(B\backslash  A)A\big)$ and  $BM  =
 G\backslash \big(B(A\backslash  B)B\big),$  
 \item [(v)]  $G=ABAB=\<A,B\>,$ and
  \item [(vi)]  $A$ and $B$ uniquely determine $M $ in $G=\<A, B\>$.
 \end{itemize} \rm\vspace{-2pt}

\proof
(i)   By \remref{ideals}(i,iii),  $A\cap B$  fixes ${\mathfrak  i}$ and ${\mathfrak  I}$, so  lies in  $ A\cap M$; and $A\cap M$ fixes the   line on $A$ and ${\mathfrak  i}$,
so lies in $  A\cap B$.  Similarly, $A\cap B=B\cap M$.

(ii)  See (i) in the proof of \thmref{main}.
  \vspace{1pt}
 
  (iii)
 \remref{ideals}(ii,iii).
  \vspace{1pt}

  (iv)     The   points   $Ag$ ($g\in G$)   are the points $Am$ ($m\in M$) of ${\mathfrak  I}$,
   together with  
the   points $Aba\ne A $ ($a\in A,b\in  B\backslash A $) on the lines $Ba$ on  $A$.%
  \vspace{1pt}
   
   (v)   If $m\in M\backslash A$ then $Amb\not\subseteq AM$ for $b\in B\backslash A$, 
   so $Amb\subset  ABA$ by (iv).
  \vspace{1pt}
   
(vi)  Use (iii,iv,v).\qed 
 
  \Theorem
  \label{converse}  Every projective plane of order $n$ having a collineation group with a flag-orbit of size $n^3$
  arises  via {\thmref{main}.} 
 $($In particular$,$  every soft plane   arises   via \thmref{main}.$)$
  
  \proof
  Given a projective plane  $\bar\pi$ of order $n$ and a collineation group $G$ with a flag-orbit of size $n^3$, we need to show that $G$ is a soft group with respect to a distinguished flag  $(\bar \infty,L_{\bar \infty})$, 
  we need subgroups $A,B$ and $M$ behaving as in \thmref{main} that produce  a plane  $\pi$, and finally 
  we must prove that  $\pi \cong \bar\pi$.  
  
  \vspace{1pt}
   (i) {\em Soft group.} 
   For a flag $(\bar A,\bar B)$  in the given flag-orbit  let $A$ and $B$ be the stabilizers  
   of $\bar A$ and $\bar B$ in $G$,
   respectively.
Then $n^3=|G\col \!A\cap B| = |G\col\! B | |B\col \!A\cap B| $ with $|G\col \!B |\le n^2+n+1$ and $ |B\col\! A\cap B| \le n+1$.
Since $|G\col\! B |$ and $ |B\col\! A\cap B| $  divide  $n^3$ it follows that 
 $|G\col \!B | = n^2$ and  $ |B\col\! A\cap B| =n$, so  $|G|=n^3k$ with $k:=|A\cap B|$.
 Similarly,    $ |A\col\! A\cap B| =n$   and  $|R|=|G\col \!A | = n^2,$ where $R:=\bar A^G$.

 Clearly $|\bar B^g\cap R|$ is the same for all $g\in G$.
 Then $r:=|\bar B\cap R|$ is $n$ or $n+1$ since 
$ \bar B\cap R$ contains the $n= |B\col\! A\cap B|$ points $\bar A^b,$ $b\in B$.
 Also,
 $\bar A$ is on the $n$  lines   $\bar B^a, $ $a\in A,$ so
 $n^2=|R|\ge |\bigcup\{\bar B^a\cap R \mid a\in A\}|=1+n(r-1)$.  Then  
  \vspace{1pt}
   $r=n$ and $\bar B$ has a unique point not in $R$.    
 Dually, $\bar A $ is on a unique  line $\bar {\mathfrak  I}\notin \bar B^G$.
    
    If $g\in G$ then $\bar A^g$ is on a unique line $ \bar {\mathfrak  I}^g$
 of $ \bar {\mathfrak  I}^G,$   so $( \bar {\mathfrak  I}\cap R)^G$ is a partition of $R$.
Since $\bar B^G\cup \bar {\mathfrak  I}^G$ is the set of lines meeting $R$ and has size 
$n^2+n^2/|\bar {\mathfrak  I}\cap R| \ne n^2+n+1, $
there is a line $L$ not meeting $R$.  Since $|R|+|L|=n^2+n+1$,  the line $L$ is unique and fixed by $G$.  Dually, $G$ fixes a point $p$, so $p\in L$.
Having a flag-orbit of size $n^3$ implies that  $G$ is a soft group  with respect to the flag~$(p,L)$.%
 
\smallskip

   (ii) {\em Behavior of $A,$ $B$ and $M$.}    We  saw that $|G|=n^3k$ and  $|A|=|B|= nk$, proving most of~\thmref{main}(1).
   \vspace{1pt}
We are given a distinguished flag  $(\bar \infty,L_{\bar \infty})$ of $\bar\pi$,  and
$G$ is flag-transitive   on the geometry $\bar\Pi$ whose points are those of $\bar\pi$ not on 
  $ L_{\bar \infty}$ and whose lines are those of $\bar\pi$  not on $\bar \infty$.
        
  If    $T:=G_{\bar {\mathfrak  I}} $, the point-transitivity of $G$ on  $\bar\Pi$
  implies that     $|G\col T|=|\bar{\mathfrak  I}^G|=n$, so   $|T|=n^3k/n$. 
Let   $\bar {\mathfrak  i}$  be the point  
\vspace{1pt}
  $\bar B \cap L_{\bar \infty}$ of $\bar\pi$   and 
 $M:=G_{\bar {\mathfrak  I} \bar {\mathfrak  i}}=T_{\bar {\mathfrak  i}}$  (compare  \remref{ideals}).
Since $B$ is transitive on $\bar B\backslash\{\bar {\mathfrak  i}\}$ it is transitive on the  $n$
lines $\ne L_{\bar {\infty}}$ on $\bar \infty$; so is $G_{\bar {\mathfrak  i}}\ge B$,
and then $G_{\bar {\mathfrak  i}}= B(G_{\bar {\mathfrak  i}})_{\bar {\mathfrak  I}} =BM $.
\vspace{1pt}
 The transitivity of  $G$, $T$  and $A$   on $L_{\bar \infty}\backslash\{\bar \infty\}$ implies that
   $|G_{\bar {\mathfrak  i}}|= |G|/|\bar{\mathfrak i\,}^G|=n^3k/n$,
 $|M|= |T|/|\bar{\mathfrak i\,}^T|=n^2k/n $, 
$G=AG_{\bar{\mathfrak  i}}=A(BM)$ and  $T=AT_{\bar{\mathfrak  i}}=AM $,
proving (1,2,3).   
   
  Let $\Pi$ denote the geometry with points and lines the cosets of $A$ and $B$ and 
   incidence determined by nonempty intersection. 
   The map 
  $\varphi\col \bar A{}^x\mapsto Ax, {\bar B{}^x\mapsto Bx } $~for $x\in G$  induces an isomorphism  
  $\bar\Pi\to \Pi$. (This is well-defined: if   $\bar A{}^x=\bar A{}^{x'}$ then 
  $x'x^{-1}\in G_{\bar A}=A$. It is an isomorphism:   if  $\bar A{}^x$ is on $\bar B{}^y$
  \vspace{1pt}
  then  some ${g\in G}$~sends~$\bar A\mapsto \bar A^x$ and $\bar B\mapsto \bar B^y$,
  so $g\in G_{\bar A}x\cap G_{\bar B}y= 
  (\bar A{}^x)^\varphi \cap (\bar B{}^h)^\varphi .$)
Using \cite[Lemma~4]{HM}~(or~the proof of \thmref{main}(iii)),   (4) holds since any two points of $\Pi$ are on at most~one~line.%
 
\smallskip

(iii)  {\em Isomorphism.}  
First note    $(*) $: By \remref{ideals}(ii),  $G_{\bar{\mathfrak  I}}=T=AM=G_{\mathfrak  I}$ and      
    $G_{\bar {\mathfrak  i}} =BM=G_{\mathfrak  i}$; 
      and  ${\bar{\mathfrak  i}}$ is on the lines  
      $\bar B^x, x\in G_{\bar {\mathfrak  i}}$, while 
      ${{\mathfrak  i}}$ is on the lines 
      $ B x, x\in G_{ {\mathfrak  i}}$, with similar statements for 
      $\bar {\mathfrak  I}$  and $\mathfrak  I$.
      (For, $G_{ \bar{\mathfrak  i}}$ is transitive on the lines of $\bar\Pi$ on  ${\bar{\mathfrak  i}}$, while $G_{ {\mathfrak  i}} =BM$ is transitive on the lines of $\Pi$
      on   ${{\mathfrak  i}}$ by \remref{ideals}(i).)

 By  \thmref{main}, the subgroups $A,B$ and $M$ of $G$  produce  a projective plane $\pi $ 
  on which  $  G$ acts as a soft group with respect to a flag $(\infty,L_{\infty})$.
By $(*)$, extending $\varphi$ by sending  
 $\bar\infty    \mapsto \infty$,
 $ L_{\bar\infty}\mapsto L_{\infty}$,
  $\bar {\mathfrak  I}{}^x\mapsto  {\mathfrak  I}x$  and
    $\bar {\mathfrak  i}{}^x\mapsto  {\mathfrak  i}x$ is well-defined and
   yields the desired isomorphism    $\bar\pi \to \pi$.
\qed    
   
 \corollary
   \label{not p-group} 
   If $n$ is a  power of a prime $p$ then a Sylow $p$-subgroup of $G$ can be used in
   {\thmref{main}} to produce $\pi$.   If $n$ is prime then $\pi$ is desarguesian.
   
\rm
\proof
Since $G$ is transitive on the $n^3$ flags opposite $(\infty, L_\infty)$ 
the same is true of a Sylow $p$-subgroup $P$ of $G$.  Now apply \thmref{converse}
to $P$ for the first statement.

If $n$ is prime then we may assume that  $G$ is an $n$-group.
If $|G|=n^3$   it is straightforward to determine both $G$ and $\pi$. 
If $|G|>n^3$ then some   $g\in A\cap B$ has order $n$ (by \remref{faithful}).
  Since $g$ fixes  $A $ and  the point on both $B$ and $L_\infty$
  it acts on $n-1$ points of $B$ and so fixes  $B$ pointwise; 
  and dually it fixes all lines $Bx$ on $A$.
Similarly, $g$ also fixes  each line $Bx$ pointwise, so $g=1$, which is not the case. \qed

\smallskip\smallskip

\remref{solvable}  contains examples of this corollary.
If $G$ is solvable then the same elementary argument shows that a 
Hall subgroup  of $G$ for the primes dividing $n$
is transitive on the flags opposite $(\infty, L_\infty)$  and 
determines $\pi$ using  \thmref{converse}.
(In particular, {\em if  $G$ is solvable   and  
$n$ is odd  then it may be assumed that   $|G|$ is odd.})

 \smallskip
    \thmref{main}(4)  is equivalent to  {\em super-noncommutativity}:

\smallskip

(4$'$) If $  a\in A \backslash B$  and  $  b\in B\backslash A$
 then  $ab\ne b'a'$ whenever $b'\in B, a'\in A$.
 \smallskip
        
        One instance of this   follows easily from the other conditions in \thmref{main}:       
  \lemma
  Assume that $(1,2,3)$ hold.
If $     a\in A   ,$ $b\in B $  and $ [a,b] \in AM\backslash A   $
 then     $ab\ne b'a'$ whenever $b'\in B,$ $ a'\in A$.   
   
\smallskip  \proof\rm
   If  $x:=[a,b] \in  AM\backslash A $ and
   $ab=b'a'$ with $b'\in B,$ $ a'\in A$, then $b'a'=bax$. Now
   $b^{-1}b'=axa'^{-1}\in AM\cap B=A\cap B$ 
   (proved exactly as in step (i) of  the proof of \thmref{main}),
   whereas $x\notin A$.\qed
 
   \lemma
  \label{B^a} 
  If $ a\in A\backslash B$  then  $B^a\cap B\le A\cap B.$
     
  \smallskip  \proof\rm
  If $b=b'^a$ for $b,b'\in B$ then $ab=b'a\in AB\cap BA=A\cup B$
   by (4), so  $ab\in A$
  since $b'a\notin B$  .  \qed

\section{Some properties of  $\pi$   and  $G$}
 \label{properties} 
This  section contains   
results  concerning a soft plane $\pi$  and   a soft group $G$   of collineations of
$\pi$  
    in Theorems~\ref{main} and  \ref{converse}.%

      \subsection{Group  summary}
 \label{summary}
 First we summarize group-theoretic results. 

     \smallskip

1.   If $B M \nor G $ then $BM$ contains a normal subgroup of $G$ that is elementary abelian   of order $n^2$
(Proposition~\ref{BM}).

2.   If $BM$ is abelian then  it is  normal in $G$ and elementary abelian  of order $n^2$ 
(Proposition~\ref{BM abelian}).

3.  
If  $AM\nor G$ and $BM\nor G$, or  if $[A,B]\leq M,$ then 
 $\pi$ is a projective plane over a semifield$,$ and $G$ has a  normal
 subgroup behaving  as in \exsref{Heisenberg}  (\corref{characterize semifields}).

4.  If $N\nor G$ with $N$ abelian and $|N|\ge n^2$  then $|N|=n^2 ,$ and 
either (i)  $G=A   N  $ with $A\cap N =1$  or (ii)
$G= B N$ with $ B\cap N=1$.  Moreover, 
  either $N$ is elementary abelian or  both (i) and (ii) hold
  (\propref{abelian normal}).

  5.   If $n$ is even, and either $k$ is odd or $n$ is not a square, then the involutions in $G$ all
  lie in the union of two normal elementary abelian 2-subgroups whose product has class 2
  (\corref{n even}).

 6.   
   If $n\equiv 2$ (mod  {\rm  4}$)$ then $n=2$
   (\propref{2 mod 4}).

7. 
If $M\nor G$ then $A\cap B=1,$  and $|G|=n^3$   is either odd or a power of $2$
  (\propref{M elations}(ii,iii)).
 
8.  
If $M\nor G$ and $C_A(M)\ne1$ then $M$ is elementary abelian.
Thus, if $M\nor G$ and $n$ is not a prime power then, by conjugation,
 both  $A$ and $B$ 
induce  groups of $|M|$ automorphisms of  $M$   (\propref{is prime power}). 

9.  If $n$ is a power of a prime $p$ then it may be assumed that $G$ is a $p$-group and $n>p$
 (\corref{not p-group}). 

10.  If $G$ is solvable and $n$ is odd then it may be assumed that  $|G|$  is odd 
 (following \corref{not p-group}). 

\medskip
 It seems likely that $G$ has to be solvable.
 \secref{appendix} contain  a little  more 
 information concerning the possible   normal structure  of $G$.

Note that  \cite{BR} contains  a significant restriction on $n$ that is not group-theoretic: 
if $n\equiv 1$ (mod~4) then $n$ is a sum of two squares  (by 6.~the possibility 
 $n\equiv 2$ (mod~4)  in \cite{BR} is not relevant).

 \subsection{Elements fixing a line pointwise {\rm (cf.~\cite[Secs.~3.1,\,4.3]{De})}}
  \label{Elations}
For a line $W$ and a point $w$ of $\pi$ let $\dot\Gamma (w,W)$ denote the group of all elements of $\Aut\, \pi$ 
 that fix all points on $W$ and all lines on $w$; these are called {\em$ (w,W)$-elations} if $w$ is on $W$ and {\em$(w,W)$-homologies} otherwise. Here $W$ is an {\em axis} and $w$ is  a {\em center }of 
 each element of $\dot\Gamma (w,W)$;  they are unique  for a nontrivial element. 
Let $\dot\Gamma (L_\infty)$  and $\dot\Gamma (\infty)$ denote the groups of  elations  in $\Aut\, \pi$ with axis 
$L_\infty$ resp. center~$ \infty$.  If 
$\dot\Gamma (L_\infty)$ is transitive on the $n^2$ points not on $ L_\infty$ then $\pi$ is  a {\em translation plane} (with  {\em translation group}  $\dot\Gamma (L_\infty)$).
By
  Propositions \ref{Heisenberg decorated}(ii)~and~\ref{likeable decorated}(ii)   
  the plane  in \thmref{main} can be a translation plane whose 
 translation group  is not contained in $G$.

Instead of $\dot\Gamma (w,W)$ we will    focus on $\Gamma (w,W) := \dot\Gamma (w,W) \cap G\hspace{-.4pt}$  and on the  normal subgroups
  $\Gamma (L_\infty):= \dot\Gamma (L_\infty) \cap G$ 
 and $\Gamma (\infty):= \dot\Gamma (\infty) \cap G$  of $G$. 
 
   \Proposition     
 \label{BM}
 If  $\<B^A\>\leq BM$ $($in particular$,$ if  $B M \nor G ),$  then  $n$ is a prime power$,$
 $B M \nor G $  and  $\pi$  is a translation plane  
 with translation group $\Gamma (L_\infty) \le BM.$%
 \rm     \smallskip

 \proof
 Since $B^A\subseteq BM$, if $a\in A$ then $(B^A)^a=B^A$ fixes ${\mathfrak  i} a$ by \remref{ideals}(ii). Then    
 $\<B^A\>$ fixes all points of $L_\infty$.  If  $\<B^A\>$ contains a nontrivial homology $h$ then
 its center is not on $L_\infty$.  By \remref{ideals}(iv) and  \cite[4.3.2]{De}, 
  $\pi$ is a translation plane with translation group 
 $\Gamma(L_\infty) < \<h^{BM}\> \le BM$.       
 
If  $\<B^A\>$ consists of elations then  all fixed points of $B$ are on $L_\infty$, so $A\cap B=1$
and $|B|=n$ by \thmref{main}(1).  
Since $B$ fixes the line $B$ it is    $\Gamma({\mathfrak  i},L_\infty)$.
 Then $\pi $ is a translation plane with translation group 
 $\<B^A\>$   \cite[3.1.20]{De},  so  $\<B^A\>=BM$ since both groups have order $n^2$ 
 by \thmref{main}(2).  \qed

\corollary
\label{characterize semifields}
If either of the following holds then 
 $\pi$ is a projective plane over a semifield %
 and $G$ has a  normal  subgroup 
behaving  as in \exsref{Heisenberg}$:$

\begin{itemize}
\item[\rm(i)] $AM$ and $BM$ are normal in $G,$ or 
\item[\rm(ii)] $[A,B]\leq M$.
\end{itemize} \rm

\proof
(i)  Use the proposition together with  \cite[(3.1.22)(f)]{De} and \cite{Hi};
the stated normal subgroup is $\<\Gamma(\infty), \Gamma(L_\infty) \>$.
  \vspace{1pt}
 
(ii) Since $B^A\subseteq BM$ and $A^B\subseteq AM$ 
this follows from (i) and the proposition.\qed%

 \proposition  
 \label{BM abelian}
  If $BM$ is abelian then 
 $n$ is a prime power and  $\pi$  is a translation plane  with translation group  $BM$. \rm
 \smallskip
 
 \proof
 By \remref{ideals}(iv), $BM$ is transitive on the $n^2$ points not on $L_\infty$, so
$|BM|=n^2$   (by \remref{faithful})  and   $k=1$ (by \thmref{main}(2)).
 Since $BM$ is transitive on the  $n$ lines $\ne L_\infty$ on $\mathfrak i$
 (by \remref{ideals}(i,ii)),
 the kernel of that action has~order $|BM|/n=n $   and   consists of elations: it is  
 $ \Gamma({\mathfrak i},L_\infty)$.    As in   the proof of \propref{BM}, $\pi $ is a translation plane with translation group 
$\< \Gamma({\mathfrak i},L_\infty)^A\>=BM$.\qed

\proposition
\label{abelian normal}
If $G$ has an abelian normal  subgroup $N$ of order $\ge n^2$ 
then  either $\pi$   or its dual is a translation plane with translation group  
$N,$  or
   $|N|=n^2,$ $N$ is transitive on  both the points not on $L_\infty$ and  the lines 
   not on $\infty  ,$ and
     $G=A   N =B N$ with $A\cap N=B\cap N=1$. \rm

\smallskip
\proof
Most of the proof consists of showing  that {\em    $\pi$ is a translation plane
with translation group  $ N=\Gamma(L_\infty)$
if  $N_B\ne1$.}~First note that 
  $N_B$ cannot fix  every point on~$B$.  For otherwise, $B$ would be an axis of every element of 
  $N_B$ and hence be fixed by the abelian group $N$; but then $N=1$ by 
  the  transitivity of $G $
  since $N\nor G$ (together with \remref{faithful}). 

 Since $G_B$ is transitive on the set of points    of  $B$ not on  $L_\infty$, that set  is the union of  nontrivial 
   $N_B$-point-orbits, each of which uniquely determines $B$.   Any $x\in N$  sends this union to the union of point-orbits
   of $(N_B)^x=N_B$, each of which  determines the line
    $Bx$ of $\pi$.  Then $Bx$  either is $B$ or  is disjoint from $B$ outside $L_\infty$, so 
    it is on the point ${\mathfrak  i}$  of  $ L_\infty$ on $B$.
   
   Since $N$ fixes   $L_\infty$ and permutes   other lines on ${\mathfrak  i}$ it fixes 
   ${\mathfrak  i}$.   
   Since   $N\nor G$
the transitivity of  $G$   on $L_\infty\backslash\{\infty\}$  implies that $N$ fixes every point of $L_\infty$.  No nontrivial element of $N$ can be a homology with axis 
   $L_\infty$ since the abelian group $N$ would have to fix its center.
    Then $N $  lies in $ \Gamma(L_\infty)$, which has order   $\le n^2\le |N|$
     \cite[4.3.2]{De},
   so $\pi$ is a translation plane,  as asserted.

Thus,  if neither $\pi$  nor its dual is a translation plane then $N_A=N_B=1.$  
Counting the number of $N$-images of $B$ yields
$n^2\ge |N|/|N_B| \ge n^2 /1$.  Then $N$ is transitive on the lines not on 
  $\infty$, and, dually, on 
the points not on   $L_\infty$.  
   \qed
 
\smallskip\smallskip
The groups in \exsref{Heisenberg} with entries from  commutative semifields  have   abelian normal subgroups  
transitive on the stated sets of points and lines.

\smallskip
If $X$ is a subgroup of $G$  let $X_0$ denote the set of elations in $X$.
   \exsref{list} with $BM\nor G$ have
 $|A_0|=$   $q $  in (ii,iv,vi),   $ 3$ in (iii) and  $2q=16$  in  (v).
   
\lemma\label{M normalizes}
{\rm(i)}
$A_0=\Gamma(\infty,{\mathfrak  I}), B_0=\Gamma({\mathfrak  i}, L_\infty) $ and
$M_0=\Gamma(\infty, L_\infty)  ,$ 

{\rm(ii)} $M_0\nor G,$ and $M_0$ 
contains every  subgroup of $M$ that is  normal in $G,$

{\rm(iii)}
$M$ normalizes $A_0$ and $B_0,$

{\rm(iv)} 
$|A_0|, |B_0| $ and $|M_0| $ divide $n,$ and

{\rm(v)} Every elation in   $G$ lies in $ \Gamma(\infty) \cup  \Gamma(L_\infty).$
\rm
 
\proof (i) $A_0$ consists of elations fixing $\infty$, $ L_\infty$ and $A$ and so  having axis ${\mathfrak  I}$ and center $\infty$. Every $(\infty,{\mathfrak  I})$-elation fixes   $A$.
 The cases $B_0$ and $M_0$ are similar.
  \vspace{1pt}
 
 (ii)     Clearly  
 $\Gamma(\infty,L_\infty)\nor G$,  and $\Gamma(\infty,L_\infty)\le 
 G_{{\mathfrak  I}{\mathfrak  i}}=M$ by \remref{ideals}(iii).  
 Let $M_1\le M <G_{\mathfrak  i}$ with $M_1\nor G$.  If  $g\in G$ then $M_1^g=M_1$ fixes  
 ${\mathfrak  i}g$, so $M_1$ fixes all points on $L_\infty$.  Dually, $M_1$ fixes all lines on $\infty$, so $M_1\le \Gamma(\infty,L_\infty)=M_0$.
  \vspace{1pt}

  (iii)   Let $a\in A_0$ and $g\in M\le G_{\mathfrak  I} $.  By (i), $a_0^g$ has axis 
 ${\mathfrak  I}g={\mathfrak  I},$  so $a_0^g\in A_0$.

  \vspace{1pt}

  (iv)  \cite[p.~187]{De}. 
  \vspace{1pt}
  
  (v)       Every  nontrivial elation   with center $w\ne \infty$ fixes $L_\infty$,   so
    $w\in L_\infty $; its   axis  contains $\infty$ and  $w$ and so is $L_\infty$.  
\qed%

    \smallskip\smallskip
 If $s$ is an integer and $p$ is a prime  then  $s_p$ denotes the largest $p$-power  dividing~$s$.%
  
\proposition\label{L0}
\label{n2}  
If $B_0\ne1$ then
\vspace{-2pt}
\begin{itemize}
\item[(i)] $\Gamma(L_\infty)$\ is an   elementary abelian  $p$-group  for some  prime~$p,$

\item[(ii)]   $\{  B_0^a \backslash\{1\}\mid  a\in A \}$ is  a $G$-invariant partition of
$\Gamma(L_\infty) \backslash  M_0$  into $n$ subsets$,$   
\item[(iii)]
$ n/|M_0|  =(|\Gamma(L_\infty) \!: \! M_0|-1)/(| B_0|-1)$  and 
$|M_0|=n_p\ge |B_0|>1 ,$

\item[(iv)]  If $n$ is not a $p$-power  then $n$    determines $|B_0|$ and 
$|\Gamma(L_\infty)|,$  and

\item[(v)]  If $n$ is not a $p$-power and $A_0\ne1$ then  
$|A_0|=|B_0|,$
$|\Gamma(\infty)|=|\Gamma(L_\infty)| ,$ \
 $\Gamma(\infty) \Gamma(L_\infty) $ 
     is a  class $2$ group whose center is
     $M_0 = { \Gamma(\infty) \cap \Gamma(L_\infty)}  ,$  and 
elements of $\Gamma(\infty)\backslash M_0$ and $\Gamma(L_\infty)\backslash M_0$ never  commute.

      \end{itemize}
\rm

  \proof
(i)  Since $B_0^a=\Gamma(\mathfrak  i a,L_\infty)$ for $a\in A$, 
there are nontrivial  elations in $\Gamma(L_\infty)$ with   different centers, so this follows from 
\cite[4.3.4(b)]{De}.
  \vspace{1pt}

 (ii)
    Every element of 
$\Gamma(L_\infty) \backslash M_0$ is an elation whose center is in 
$L_ \infty\backslash\{\infty\}$ and so is 
$\mathfrak i a$ for some $a\in A$.  
Then  $\{  B_0^a \backslash\{1\}\mid  a\in A \}$ is  a  partition of
$\Gamma(L_\infty) \backslash M_0 $ of size $|{\mathfrak i}A|=n$.  
 It is $G$-invariant since $  G$ acts on $L_ \infty\backslash\{\infty\}$. 
  \vspace{1pt}
   
(iii)  By (ii), 
$|\Gamma(L_\infty) |-|M_0|\! =\! n(|B_0|-1)$, so 
  $  n/|M_0|  \! = \!  \big( |\Gamma(L_\infty)\!: \!M_0|-1\big  )  /(|B_0|-1)$.
  By   (i) and 
  \lemref{M normalizes}(iv),  $  n/|M_0| $ is an integer not divisible by $p$.
  Then $|M_0| =n_p$ since $|M_0|$ is a  $p$-power, 
  and  $n_p\ge |B_0|>1$, again by \lemref{M normalizes}(iv).
  \vspace{1pt}

(iv)   
By (iii),  $ 0< (n/|M_0|) -1   =  
 |B_0|\big((|\Gamma(L_\infty)|/|M_0||B_0|  )  -   1   \big)   /(  |B_0| -1  )$.
Since $\Gamma(L_\infty)\ge M_0B_0\cong M_0\oplus B_0$,   
 \  $|\Gamma(L_\infty)|/|M_0||B_0| $
 is a $p$-power  and is not 1.  Then  $|B_0|$ is the largest $p$-power dividing 
 $(n/|M_0|) -1 \!= \!(n/n_p) -1$,  which  determines~$|\Gamma(L_\infty)|$.%
  \vspace{1pt}
  
(v) 
The order equalities follow from (iv).
By (i),  
$ \<\Gamma(\infty), \Gamma(L_\infty) \>=  \Gamma(\infty)  \Gamma(L_\infty) $ 
is a $p$-group with 
 $ \<\Gamma(\infty), \! \Gamma(L_\infty) \>' \!  = \! [ \Gamma(\infty), \Gamma(L_\infty)]$ $\le \Gamma(\infty)\cap \Gamma(L_\infty)
\!  \le\!  Z \big ( \Gamma(\infty)  \Gamma(L_\infty) \big),$ 
so   $ \Gamma(\infty) \Gamma(L_\infty) $  has class 2. 
Also,  $M_0 =  \Gamma(\infty, L_\infty) 
 \le Z \big ( \Gamma(\infty)  \Gamma(L_\infty) \big). $ 

If $x\in \Gamma(\infty)\backslash M_0$ and $y\in \Gamma(L_\infty)\backslash M_0$ 
commute then $y$ fixes the axis of $x$,
 whereas $L_\infty$ is the only line on $\infty$  fixed by $y$.  
Consequently,   $ M_0=Z\big( \Gamma(\infty) \Gamma(L_\infty) \big)$.  
 \qed  
 
 \corollary
   \label{n even} 
   If $n$ is even$,$ and either $k$ is odd or $n$ is not a square$,$ then all involutions of $G$ lie in $\Gamma( L_\infty)\cup\Gamma(\infty),$ where
   $\Gamma(L_\infty)\ge B_0\times M_0$ and $\Gamma(\infty)\ge A_0\times M_0$
   are  elementary abelian   $2$-groups and $\Gamma(\infty)\Gamma(L_\infty) $ has class $2$.
\rm

\proof
If $k$ is odd or $n$ is not a square then  all involutions are elations \cite[4.1.9]{De} and hence lie in 
   $\Gamma(\infty)\cup\Gamma(L_\infty) $  (by \lemref{M normalizes}(v)).
The remaining assertions follow from  \lemref{M normalizes}(i)  and  \propref{L0}(i,v).\qed

 \proposition
   \label{2 mod 4} 
   If $n\equiv 2$ $(mod $ {\rm  4}$)$ then $n=2$.

   \rm
   \proof   Since  $n$ is not a square, $|A\cap B|$ is odd    \cite[4.1.9]{De}. Then
  ${|G|_2=n_2^3=8}$~by \thmref{main}(1).  
   By \corref{n even},
 $D:= \Gamma( \infty) \Gamma(L_\infty) $ is a normal dihedral subgroup of $G$.
\!\!\! Any $a\in A$ of odd order centralizes $D$, where
  $D\ge \Gamma(L_\infty)\ge B_0$ and $A\cap B_0=1,$   so  $a\in A\cap B$ by  
 \thmref{main}(4).  Then  $|A|/2$ divides $|A\cap B|=|A|/n$, so  $n=2$.\qed%
 
 \subsection{\!The case $M\nor G$.}
  \label{normal}
{\em \! We  now assume  that the subgroup $M\!$ occurring in 
 {\rm Theorem \ref{main}} is normal  in $G$.}
 All of the planes $\pi$  in \exsref{list} arise from   soft groups 
  $G\le \Aut\, \pi$ for which  $M$  behaves in this manner, but in general $M$ need not be normal in $G$  (Propositions~\ref{Heisenberg decorated}(i) and~\ref{likeable decorated}(i)).

 \Proposition
\label{M elations} 
{\rm(i)}
$M=M_0=\Gamma(\infty,L_\infty),$    \vspace{-2pt}
\begin{itemize}
\item[\rm(ii)] $A\cap B=1,$  $   |\Gamma(\infty,L_\infty)|=   |M|=n $  and  $|G|=n^3,$

\item[\rm(iii)]
$|G|$ is either odd or a power of $2,$
 and
\item[\rm(iv)]
Every element of $G$ can be written $amb$ for unique $a\in A, m \in  M, b\in B$. 
\end{itemize}
 \rm

\proof
(i)   \lemref{M normalizes}(i,ii).
  \vspace{1pt}
 
 (ii)  By \remref{ideals}(i)  and (i),  $M$   is regular on  ${\mathfrak  I}\backslash\{\infty\}$, so    \thmref{main}(1) implies that  
 $n =|M|=n k$  and~${|G|=n^3}$.%
  \vspace{1pt}

(iii)  Assume that $n$ is even.  Since $A\cap B=1,$ each involution in $B$    is an elation. 
By \propref{L0}(i), $\Gamma(L_\infty) >M $ is elementary abelian.
Then $|M|=n$ and~$|G|=n^3$ are  powers of 2.
  \vspace{1pt}
  
(iv)  Use  \thmref{main}(1,2,3) and (ii). \qed
   
       \proposition
     \label{is prime power}
  If  $C_B(M)\ne1$ then $M$ is elementary abelian of order $n$.

If  $n$ is not a prime power then$,$
by conjugation$,$
 both  $A$ and $B$ 
induce  groups of $|M|$ automorphisms of  $M$.  
 \rm\smallskip
  
  \proof 
Let $ 1\ne c\in C_B(M)$.   Then $c$ fixes each line on    $\mathfrak  i$     
 since $(Bm) c=Bcm=Bm $ for $m\in M$ (cf. \remref{ideals}(i)). 
  Homologies have order dividing $n-1$,  but $|G|=n^3$ by  \propref{M elations}(ii).  Then
   $c$ is an elation,  so $c\in B_0$ and
    $M$ is elementary abelian of order $n$  by  Propositions
\ref{M elations}  and  \ref{L0}(i).
 
 If  $n$ is not a prime power then  $C_B(M)=1$, so
 $B$  induces a group of $|B|=|M|$ automorphisms of $M$.
     \qed%

 \section{Additional groups}

  \label{Additional groups}
    \label{appendix}
Examples~\ref{Heisenberg} and \ref{list}
  focussed on planes not on the groups that produce them.  
Here we deal with
 additional groups that  occur as  soft groups of planes in some of those examples.
For $G$ in  Examples~\ref{Heisenberg} 
or \ref{list}(ii) (using   groups in  \cite{Ka}
and provided below in the  proof   of \propref{likeable decorated}), 
we assume that 
 the underlying field or semifield has a nontrivial automorphism $\a$    that acts on matrix entries
 (for example, see the start of the proof of \propref{likeable decorated}).
 This induces a collineation (also called~$\a$) of the associated plane $\pi$.
 
\Remark
\label{solvable}\rm
The group $\tilde G:=G\<\a\>$ acts on  $ \pi$ as a soft group 
(using Examples~\ref{Heisenberg} or \ref{list}(ii) together with  \eqref{table}).  
 If $ \<\a\>$ is  not a $p$-group then  $G\<\a\>$  also is not a $p$-group:
 {\em  $G$ can be solvable but not nilpotent in  {\thmref{main}}}.
 \smallskip\smallskip
   
 The stabilizer $\tilde G_\Phi$ of the flag $\Phi:=(A , B) $ of $\pi$ is  $\<\a\>$.
The nonabelian groups  $\tilde A:=A\<\a\> , $   $\tilde B:=B\<\a\>  $ and 
$\tilde M:=M\<\a\>  $
 play the roles of $A,$   $B$ and $M$ for~$\tilde G $.
 
  In the rest of this section {\em we assume that $\a$ has order $p$}.  
Then \exsref{list}(ii,iv) with $q=2^{2e+1}$ cannot occur.

 \proposition
\label{Heisenberg decorated}
Given any prime $p,$
for infinitely many  $p$-powers  $n$ there are~both desarguesian planes and nondesarguesian semifield planes having   collineation groups  behaving as follows$:$
 \begin{itemize}
 \item[\rm(i)]  $\tilde{G}\le \Aut\, \pi $ is a soft $p$-group  
 of order $n^3p,$  with  corresponding subgroups   
    $\tilde{A} \hspace{-.5pt}\tilde{M} ,$       $\tilde{B} \hspace{-.5pt} \tilde{M}
    \ntrianglelefteq\tilde{G}$   
  and $\tilde{M}\ntrianglelefteq \tilde{A} \tilde{M},$   $\tilde{B} \tilde{M};$
  and  
   
 \vspace{2pt}
 \item[\rm(ii)]   $\ddot{G}\le \Aut\, \pi $ is a soft $p$-group  of order $n^3,$ not isomorphic to  any group  in {\rm\exsref{Heisenberg}}$,$   
with corresponding subgroups   
  $\ddot{A} \ddot{M} ,\ddot{B}\hspace{-.5pt}\hspace{-.5pt} \ddot{M} \ntrianglelefteq\ddot{G}$  while $\ddot{M}\nor \ddot{G}$.
 $($Moreover$,$ the translation group of the associated translation plane is not 
 contained in 
 $\ddot{G}.)$
  \end{itemize}
 \rm 
   
 \proposition
\label{likeable decorated}
For any given odd prime $p\equiv 2$ $($mod $3),$
for infinitely many  $p$-powers  $n$ there are planes  $\pi$  in {\rm\exsref{list}(ii)}
 having   collineation groups  behaving as follows$:$
 \begin{itemize}
 \item[\rm(i)]  $\tilde{G} \le \Aut\, \pi $ is a soft $p$-group  
 of order $q^6p,$   with  corresponding subgroup    
   $\tilde{B} \hspace{-.5pt} \tilde{M}\ntrianglelefteq\tilde{G};$   and
  
   \vspace{1pt}
 \item[\rm(ii)]   $\ddot{G}\le \Aut\, \pi $ is a soft $p$-group  of order $q^6,$ not 
 isomorphic to any group   in {\rm\exsref{list}(ii)}$,$
 with corresponding subgroups   $\ddot{A} \ddot{M},$       
 $\ddot{B}\hspace{-.5pt}\hspace{-.5pt} \ddot{M}
 \ntrianglelefteq\ddot{G}$  and $\ddot{M}\nor \ddot{G}$.
 $($Moreover$,$ the translation group of the associated translation plane is not 
 contained in  $\ddot{G}.)$
  \end{itemize}
 \rm

  \smallskip   
  
  Since the proofs of these propositions are very similar  we will only sketch
  a proof of the second one.
  \smallskip   
 
 \proof
Let $F:=\F_q$ and $G:=  A\ltimes F^4 ,$ where the elementary abelian group 
  $A:=\{ A(t,u)\mid t,u\in F\}$ acts on $F^4 $  via 
        \begin{equation} \label{commutator}
     [A(t,u) , (x,y,z,w) ] =-\big (0, xt, xu+yt,xf(t,u)+yu+zt \big)
     \end{equation}
    with $f(t,u):=tu-\frac{1}{3} t^3+l(t)$  \vspace{1pt}
    for a suitable additive map $l\col F \to F$ (cf. \cite{Ka}). (Here $l=0$ if  the 
    final  characteristic 5 instances of  \exsref{list}(ii) are  ignored.)
      \exsref{list}(ii) are obtained using  $G$ 
       and  the subgroups  $A$, $B:=(F,F,0,0)$ and $M:=(0,0,F,F)$.
The  action of $\a$ on $G$ is $A(t,u)(x,y,z,w)\mapsto A(t^\a,u^\a)(x^\a,y^\a,z^\a,w^\a)$.

\smallskip
(i)
We   noted that 
   $\tilde G=G\<\a\>$  is a soft $p$-group.   
   If $a^\a\ne a\in A$ then $\a^a=\a[\a,a]\notin BM\<\a\>$, so
    $\tilde M=M\<\a\>$  and  
    $\tilde B\tilde M=BM\<\a\>$ are not normal in $\tilde G$.
  \vspace{1pt}
 
 (ii)    By    \eqref{commutator},  $G'=[A,F^4]=(0,F,F,F)$. Then  
    $\tilde G= AG'(F,0,0,0)\<\a\> $   implies that $\tilde G/G'\cong A(F,0,0,0)\<\a\>$,
 where  $A(F,0,0,0)$ is elementary abelian by \eqref{commutator}.  
Here  $\tilde G \ron AG' $  and    
  $\tilde G/AG'  \cong  (F,0,0,0)\<\a\>$.
  
  Let $H$ be an $\a$-invariant subgroup of index $p$ in $(F,0,0,0)$,
  so $\tilde G/AG'H\cong ((F,0,0,0)/H)\times \<\a\>\cong \Z_p^2$.
   Let ${\ddot{G}} $   be a  subgroup  of index $p$ in  $\tilde G$ 
   containing $AG'H$  and such that ${\ddot{G}}\ne AG'(F,0,0,0),$ $AG'H\<\a\>$.

  We will show that $\ddot{G}$ behaves as required.

First note that  $\ddot{G}$ is a soft group for  $\pi$:
  $|\Phi^{\ddot{G}}|=|{\ddot{G}}|/|{\ddot{G}}_ \Phi|$ with 
$\a\notin {\ddot{G}}_\Phi
\le  \tilde G_\Phi= \<\a\>$  and $|\a|=p$,
so $\tilde G_\Phi=1$ and
$|\Phi^{\ddot{G}}| =|{\ddot{G}}| =q^6,$ as claimed. \vspace{1pt}

There are subgroups  \vspace{1pt} $\ddot{A},\ddot{B} , \ddot{M}$  of $\ddot{G}$
behaving as in \thmref{converse}.
 By \propref{BM},   \vspace{1pt}
      $\ddot{B}\ddot{M}\ntrianglelefteq \ddot G$
      since   the subgroup $(F,0,0,0)$ of  the translation group  $BM$  of $\pi$ is not in $ \ddot G$.
     Also by \propref{BM}, $\ddot{A}\ddot{M}\ntrianglelefteq \ddot G$
      since $\pi$ is not the dual of a  translation plane. 
      \vspace{2pt}
      
Clearly  $  \ddot{M}=\tilde M \cap \ddot{G}=M\<\a\>\cap \ddot{G}=M\nor \tilde G$.    

Finally,  we claim that  $\ddot{G}$ is not isomorphic to   $G^\star$,
an arbitrary group  in \exsref{list}(ii).  
       By \eqref{commutator}, $C_{\tilde G}(G')=C_{\tilde G}\big ((0,F,F,F) \big)=F^4$~and
    $C_{\tilde G}(AG')=(0,0,0,F)$.
    Since $\ddot{G} $ contains an element of $G\a$   it follows that
  $C_{\tilde G}(\ddot{G} )\le (0,0,0, C_F(\a))$,
       so $Z(\ddot{G} ) $ is smaller than $Z(G^\star)=(0,0,0,F)$
       and  $\ddot{G}  \not \cong G^\star$.\qed      
        
        \section{Concluding  comments}
        This paper contains a possible method for   constructing finite projective planes 
        not of prime power order, but no such plane has been found.         
Restrictions on the group are in \secref{summary}.
  For planes of prime power order
        every known soft $p$-group has nilpotence class at most 5, but 
        it seems likely that there are soft $p$-groups  of large class.

Translation planes are coordinatized using the translation group  
and its   large point-orbit.   Soft-planes are in some sense ``coordinatized'' using a soft group and its very large flag-orbit.

{ 
}

Groups of order $q^6$ implicitly occur in \exsref{list}(ii,iv,v,vi) 
(which were what \cite{Ka,FJW,JW,JJ,BM,BJJM} were accidentally
studying).  It is not  clear what standard group-theoretic  properties those groups have
beyond   ones  seen here.

     \end{document}